\documentclass[reqno,centertags, 12pt]{amsart}
\usepackage{amsmath,amsthm,amscd,amssymb,mathtools}
\usepackage{latexsym}
\usepackage{graphicx}
\usepackage[mathscr]{eucal}
\usepackage[labelformat=empty]{caption}

\newcommand{\bbC}{{\mathbb{C}}}
\newcommand{\bbD}{{\mathbb{D}}}

\newcommand{\bbR}{{\mathbb{R}}}

\newcommand{\bdone}{{\boldsymbol{1}}}

\newcommand{\bi}{\bibitem}

\newcommand{\beq}{\begin{equation}}
\newcommand{\eeq}{\end{equation}}
\newcommand{\ba}{\begin{align}}
\newcommand{\ea}{\end{align}}

\newcounter{smalllist}

\newcommand{\comm}[1]{}

\allowdisplaybreaks

\newtheorem*{p2.1}{Proposition 2.1}

\theoremstyle{definition}

\newcommand{\overbar}{\overline}

\begin{document}

\title[Spectral Theory]{Spectral Theory Sum Rules, Meromorphic Herglotz Functions and Large Deviations}

\author[Barry~Simon]{Barry Simon$^{1}$}

\


\thanks{$^1$ Barry Simon is I.B.M. Professor of Mathematics and Theoretical Physics, Emeritus, at Caltech. His email address is bsimon@caltech.edu.}

\maketitle

Almost exactly forty years ago, Kruskall and collaborators revolutionized significant parts of applied mathematics by discovering ra emarkable structure in the KdV equation.  Their main discovery is that KdV is completely integrable, with the resulting infinite number of conservation laws.  Deeper aspects concern the connection to the $1D$ Schr\"{o}dinger equation
\begin{equation}\label{1}
  -\frac{d^2}{dx^2}+V(x)
\end{equation}
where the potential, $V$, is actually fixed time data for KdV.

In particular, the conserved quantities, which are integrals of polynomials in $V$ and its derivatives, can also be expressed in terms of spectral data (eigenvalues and scattering information).  Thus one gets a \emph{sum rule}, an equality between coefficient data on one side and spectral data on the other side.  The most celebrated KdV sum rule is that of Gardner et al.:
\begin{equation}\label{2}
  \frac{1}{\pi} \int_{0}^{\infty} \log|t(E)|^{-1} E^{1/2} \, dE + \frac{2}{3}\sum_{n} |E_n|^{3/2} = \frac{1}{8} \int_{-\infty}^{\infty} V(x)^2 \, dx
\end{equation}
where $\{E_n\}$ are the negative eigenvalues and $t(E)$ the scattering theory transmission coefficient.  We note that in this sum rule all terms are positive.

While these are well known, what is not so well known is that there are much earlier spectral theory sum rules, which, depending on your point of view, go back to 1915, 1920, or 1936.  They go under the rubric Szeg\H{o}'s Theorem, which expressed in terms of Toeplitz determinants goes back to 1915.  In 1920, Szeg\H{o} realized a reformulation in terms of norms of orthogonal polynomials on the unit circle (OPUC), but it was Verblunsky in 1936 who first proved the theorem for general measures on $\partial\bbD$  in $\bbC$  and expressed it as a sum rule.

To explain the sum rule, given a probability measure, $\mu$, on $\partial\bbD$ which is non--trivial (i.e.\ not supported on a finite set of points), let $\{\Phi_n(z)\}_{n=0}^\infty$ be the monic orthogonal polynomials for $\mu$.  They obey a recursion relation
\begin{equation}\label{3}
  \Phi_{n+1}(z)=z\Phi_n(z) - \overbar{\alpha}_n \Phi^*_n(z), \quad \Phi_0 \equiv \bdone, \quad \Phi^*_n(z) = z^n \overbar{\Phi_n\left(\frac{1}{\bar{z}}\right)}
\end{equation}
where $\{\alpha_n\}_{n=0}^\infty$ are a sequence of numbers, called Verblunsky coefficients, in $\bbD$.  Verblunsky's Theorem proves a 1-1 correspondence between such $\mu$s and such sequences.

The Szeg\H{o}--Verblunsky sum rule says that if
\begin{equation}\label{4}
  d\mu(\theta) = w(\theta)\frac{d\theta}{2\pi}+d\mu_s
\end{equation}
then
\begin{equation}\label{5}
  \int \log(w(\theta))\,\frac{d\theta}{2\pi} = -\sum_{n=0}^{\infty} \log(1-\alpha_n|^2)
\end{equation}
In particular, the condition that both sides are finite at the same time implies
\begin{equation}\label{6}
  \sum_{j=0}^{\infty} |\alpha_j|^2 < \infty \iff \int \log(w(\theta)) \frac{d\theta}{2\pi} > -\infty
\end{equation}

Simon \cite{SzThm} calls a result like \eqref{6} that is an equivalence between coefficient data and measure theoretic data a \emph{spectral theory gem}.   The one above has a spectacular (albeit perhaps hidden!) consequence.  It is known that if $\{\alpha_n\}_{n=0}^\infty  \in \ell^1$, then the measure $\mu$ is purely a.c.  In \eqref{6}, $d\mu_s$ is arbitrary, so, as long as the a.c.\ part yields a finite integral, we have that $\{\alpha_n\}_{n=0}^\infty \in \ell^2$ no matter what $\mu_s$ is.  Thus we have arbitrarily singular mixed spectrum with $\ell^2$ decay of the recursion coefficients.

In 2000, Killip and I \cite{KS} (I know the publication date was 2003, not 2000 but its the Annals -- what do you expect!) found an analog of the Szeg\H{o}--Verblunsky sum rule for orthogonal polynomials on the real line (OPRL).  One now has non--trivial probability measures on $\bbR$ and $\{p_n\}_{n=0}^\infty$ are orthonormal polynomials whose recursion relation is
\begin{equation}\label{7}
  xp_n(x) = a_{n+1}p_{n+1}(x) + b_{n+1}p_n(x) + a_n p_{n-1}(x), \qquad p_{-1} \equiv 0
\end{equation}
where the Jacobi parameters obey  $b_n \in \mathbb{R}$, $a_n \geq 0$.  There is now a bijection of non--trivial probability measures of compact support on $\bbR$ and uniformly bounded sets of Jacobi paramters (Favard's Theorem).

If
\begin{equation}\label{8}
  \ d\mu(x) = w(x) dx +d\mu_s
\end{equation}
then the gem of Killip--Simon says that
\begin{equation}\label{9}
\begin{split}
  \sum_{n=1}^{\infty} (a_n-1)^2 &+ b_n^2 < \infty \\
&\textrm{if and only if} \\
 \textrm{ess supp }(d\mu) = [-2,2], \,
   Q(\mu) &< \infty   \textrm{ and } \sum_{m} (|E_m|-2)^{3/2} < \infty
\end{split}
\end{equation}
where
\begin{equation}\label{10}
  Q(\mu) = -\frac{1}{4\pi}\int_{-2}^{2} \log\left(\frac{\sqrt{4-x^2}}{2\pi w(x)}\right) \sqrt{4-x^2} \, dx
\end{equation}
The sum rule is
\begin{equation} \label{11}
  Q(\mu) + \sum_{ \mu\left(\{E_n\} \right)>0,  \ |E_n|>2} F(E_n) = \sum_{n=1}^{\infty} \left[\tfrac{1}{4}b_n^2 + \tfrac{1}{2} G(a_n) \right]
\end{equation}
where
\begin{align}
  F(\beta+\beta^{-1}) &= \tfrac{1}{4}[\beta^2+\beta^{-2}-\log(\beta^4)], \qquad \beta \in \bbR\setminus [-1,1]  \label{12} \\
  G(a) &= a^2 - 1 - \log(a^2) \label{13}
\end{align}

The gem comes from $G(a) >0 $ on $(0,\infty)\setminus \{1\}$, $G(a) = 2(a-1)^2+\textrm{O}((a-1)^3), F(E)>0$ on $\bbR\setminus [-2,2], F(E) = \tfrac{2}{3} (|E|-2)^{3/2} + \textrm{O}((|E|-2)^{5/2})$.  To get gems from the sum rule without worrying about cancellation of infinities, it is critical that all the terms are positive.

It was mysterious why there was any positive combination and if there was any meaning to the functions $G$ and $F$ which popped out of calculation and combination.  Moreover the weight $(4-x^2)^{1/2}$ was mysterious.  Prior work had something called the Szeg\H{o} condition with the weight $(4-x^2)^{-1/2}$, which is natural since under $x=2\cos \theta$ one finds that $(4-x^2)^{-1/2} \, dx$ goes to $d\theta$, up to a constant.

This situation remained for almost 15 years during which period there was considerable followup work but no really different alternate proof of the Killip--Simon result.  This situation changed dramatically in the summer of 2014 when Gamboa, Nagel, and Rouault \cite{GNR1} found a probabilistic approach using the theory of large deviations from probability theory.

Their approach shed light on all the mysteries.  The measure $(4-x^2)^{1/2}\,dx$ is just (up to scaling and normalization) the celebrated Wigner semi--circle law.  The function $G$ of \eqref{13} is just the rate function for averages of sums of independent exponential random variables, as one can compute from Cram\'{e}r's Theorem. The function F of \eqref{12} is just the logarithmic potential in a quadratic external field which occurs in numerous places in the theory of random matrices.

In the first half of my lecture, I'll discuss sum rules via meromorphic Herglotz functions and in the second half the large deviations approach of Gamboa, Nagel and Rouault.



\end{document}